%% file: ban_gut14.tex
\documentclass[12pt]{amsart}

\usepackage{epsfig,color}

\usepackage{amssymb}

\newtheorem{thm}{Theorem}[section]
\newtheorem{prop}[thm]{Proposition}
\newtheorem{defin}[thm]{Definition}

\newtheorem{corol}[thm]{Corollary}

\newtheorem{lem}[thm]{Lemma}
\newtheorem{conj}[thm]{Conjecture}
\newtheorem{rem}[thm]{Remark}

\newtheorem{exa}[thm]{Example}

\newcommand{\mm}{{\mathcal M}}

\newcommand{\N}{{\mathbb N}}

\newcommand{\R}{{\mathbb R}}
\newcommand{\RR}{{\mathbb R}^2}

\newcommand{\Z}{{\mathbb Z}}
\newcommand{\ZZ}{\mathbb Z^2}

\newcommand{\bo}{\partial} 

\newcommand{\inn}{\text{in}}        

\newcommand{\len}{\text{length}}     

\newcommand{\out}{\text{out}}   

\newcommand{\stab}{\text{Stab}}

\newcommand{\tc}{\tilde{c}}

\newcommand{\dc}{\dot{c}}
\newcommand{\hc}{\hat{c}}
\newcommand{\dhc}{\dot{\hat{c}}}

\newcommand{\al}{\alpha}

\newcommand{\ve}{\varepsilon}

\def\qed{\hfill\vrule height 5pt width 5 pt depth 0pt}

\date{\today}

\begin{document}

\title[Blocking and Flatness]
{Secure two-dimensional tori are flat}

\author{Victor Bangert and Eugene Gutkin}

\address{Mathematisches Institut\\
Albert-Ludwigs-Universit\"at\\
Eckerstrasse 1\\
79104 Freiburg im Breisgau\\
Germany}

\email{bangert@email.mathematik.uni-freiburg.de}

\address{UMK and IMPAN\\
Chopina 12/18\\
Torun\\
Poland}

\email{gutkin@mat.uni.torun.pl,\ gutkin@impan.pl}


\maketitle

\date{\today}

\begin{center}
\begin{minipage}{120mm}

\vskip 0.1in
\baselineskip 0.2in

{{\bf Abstract.} A riemannian manifold $M$ is secure if
the geodesics between any pair of points in $M$ can be 
blocked by a finite number of point obstacles.  Compact, flat manifolds
are secure. A standing conjecture says that
these are the only secure, compact riemannian manifolds. The conjecture claims, in particular, that 
a riemannian torus of any dimension is secure if
and only if it is flat. We prove this for two-dimensional tori.

}

\end{minipage}
\end{center}

\section{Introduction}  \label{intro}
We begin by describing our setting and establishing the terminology. 
By a riemannian manifold $(M,g)$ we will always mean a 
complete, connected,
infinitely differentiable riemannian manifold. We will view geodesics in
$(M,g)$ as curves, $c:I\to M,\,I\subset\R$, parameterized by arclength.
If $I=[a,b]$, we say that $x=c(a),y=c(b)$ are the endpoints of $c$. 
If $z\in M$ is an interior point of
$c$, we say that $c$ {\em passes through $z$}.

For any pair $x,y\in M$ (including $y=x$) let $G(x,y)$ be the set of 
geodesics in $(M,g)$ with endpoints $x,y$. We say that $G(x,y)$ consists 
of the geodesics {\em joining} $x$ with $y$. A finite set $B\subset M\setminus\{x,y\}$ 
is a {\em blocking set} for $x,y$ if every geodesic in $G(x,y)$ passes
through a point $b\in B$. We will also say that $B$ {\em blocks $x$ and $y$
away from each other}.

A pair $x,y\in M$ is {\em insecure} if the points $x,y$ cannot be blocked
away from each other. A riemannian manifold is insecure 
if we cannot block $x$ away from $y$ for some points $x,y\in M$.
Thus, $(M,g)$ is {\em secure}  if any point in it can be blocked away 
from any point, including itself. See \cite{Gut05} for an explanation of 
the term ``security'' and related terminology.

Which compact riemannian manifolds are secure? The only examples so far
are the flat manifolds \cite{GS06}. Researchers in the subject
believe in the following statement \cite{BG08,LS07}.
\begin{conj}   \label{main_conj}
A compact riemannian manifold is secure if and only if it is flat.
\end{conj}
Slightly restricting the setting, we state a counterpart of Conjecture~\ref{main_conj} for tori.
\begin{conj}   \label{tori_conj}
A riemannian torus is secure if and only if it is flat.
\end{conj}
In this note we establish Conjecture~\ref{tori_conj} for
two-dimensional tori.

\begin{thm}   \label{main_thm}
A two-dimensional riemannian torus is secure
if and only if it is flat.
\end{thm}

\vspace{3mm}

Conjecture~\ref{main_conj}
holds for locally symmetric spaces \cite{GS06}. Moreover, let $(M,g)$ be a 
compact locally symmetric space of noncompact type. Then no points $x,y\in M$ can 
be blocked away from each other \cite{GS06}. This is true, in particular, for 
a compact surface of constant negative curvature. There is a direct
geometric argument that shows this. It is outlined on
page 193 in \cite{GS06}. The proof of our Proposition~\ref{no_block_prop} 
uses a similar idea. We point out, however, that the two approaches differ considerably in detail.
The discussion in \cite{GS06} crucially uses the hyperbolicity of the geodesic flow 
for a surface of constant negative curvature.

Our proof of Theorem~\ref{main_thm} uses the fact that there exists a free homotopy class of closed curves
such that the periodic geodesics in this class do not foliate the entire non-flat two-torus.
Thus, there are cylindrical regions free of geodesics in this homotopy class. 
We pick such a cylinder and let $p,q$ be an arbitrary pair of points in it.
We choose an infinite strip in the universal covering of our torus, projecting  
onto the cylinder in question.  
We construct an infinite sequence of minimal godesics 
in the strip; they project into a sequence of geodesics in the cylinder connecting $p$ and $q$.
We analyze the asymptotic behavior of minimal geodesics in the strip, as they become longer and longer. 
We prove that long minimal geodesics spend almost all of their time in the vicinity of a
boundary component of the strip.  See sections~\ref{constr} and~\ref{key}, in particular, Lemma~\ref{key_lem}.
This allows us to conclude 
that any point in the torus can block at most a finite number of geodesics in our infinite sequence
of connecting geodesics. See section~\ref{proof}. 
 
\medskip

A riemannian
manifold is {\em uniformly secure} if there exists a positive integer $s$ such that 
any points $x,y\in M$ can be blocked away from each other by at most $s$ blocking points.
The minimal such $s$ is the {\em  security threshold} of $(M,g)$.
Flat, compact manifolds are uniformly secure. Moreover, their security thresholds
are bounded above in terms of the dimension of the manifold \cite{GS06}. 
The fundamental group of a compact, uniformly secure manifold is  virtually nilpotent, and 
its topological entropy vanishes \cite{BG08}. If, in addition, the manifold has no
conjugate points, then it is flat \cite{BG08}. On the other hand,
if a manifold has positive topological
entropy and no conjugate points, then no pairs $x,y$ can 
be blocked away from each other \cite{LS07,BG08}. The crucial idea used in the
proofs of these statements is to relate the uniform blocking in $(M,g)$
with the growth, as $T\to\infty$, of the number $n_T(x,y)$ of geodesics in $(M,g)$ joining $x$ 
with $y$ and having length $\le T$. Relationships between $n_T(x,y)$ and the growth of  
$\pi_1(M)$, as well as between $n_T(x,y)$ and the topological  entropy of $(M,g)$ 
are well known \cite{Mane,Manning}.

\medskip

\noindent{\bf Acknowledgements.} A substantial portion of the present work took place during 
the  second author's visit to the Institute for Mathematical Research at ETH in Zurich.
It is a pleasure to thank the Institute and its director, Marc Burger, for hospitality and
financial support.


\section{Rays, corays, and Busemann functions} \label{general}
In this section we recall a few well known and less known facts
about rays, corays and the Busemann functions in complete,
connected riemannian manifolds of arbitrary dimensions. We will use
the notation $(M^n,g)$ for riemannian manifolds, suppressing $g$ or $n$
whenever this causes no confusion. We denote by $d(\cdot,\cdot)$
the riemannian distance on $M$. We will view geodesics as
parameterized curves $c(t),\,t\in I$, where $I\subset\R$ is a
nontrivial, possibly infinite interval, and $t$ is an arclength
parameter. We will call the set $c(I)\subset M$ the {\em trace of
$c$}.

\medskip
\begin{defin}   \label{ray_def}
\vspace{1mm}
Let $I\subset\R$ be an interval and let $c:I\rightarrow M$ be a
geodesic.
\vspace{1mm}
\begin{itemize}
\item[(a)]
The geodesic $c$ is  {\em minimal} if $d(c(t), c(s))=|t-s|$ for
all $s,t\in I$.
\item[(b)]
A {\em ray} is a minimal geodesic $c:\R_+\rightarrow M$.
\item[(c)]
Let $c:\R_+\rightarrow M$ be a ray, and let $C\subset M$ be its trace.

\noindent A ray $\tilde{c}$ is a {\em coray} to $c$ if there exists a
sequence of minimal geodesics $c_n:[0, L_n]\rightarrow M$ with
$\lim_{n\rightarrow\infty} L_n=\infty$, such that
$\lim_{n\rightarrow\infty} \dot{c}_n(0)=\dot{\tilde{c}}(0)$ and
$c_n(L_n)\in C$ for all $n\in\N$.
\end{itemize}
\end{defin}

Taking limits of minimal geodesics of finite length, we obtain the
following basic fact.

\medskip\noindent
\begin{prop}  \label{ray_lem}
A complete riemannian manifold $(M,g)$ carries a ray if and only
if it is not compact. If $c$ is a ray in $(M,g)$ and $p\in M$,
then there exists a coray $\tilde{c}$ to $c$ with
$\tilde{c}(0)=p$.
\end{prop}
\medskip
%
\begin{defin}   \label{busem_def}
Let $c:\R_+\to M$ be a ray. Its {\em Busemann function},
$B_c:M\to\R$, is defined by
\begin{equation}   \label{busem_eq}
B_c(p) = \lim_{t\rightarrow\infty}\left[d(p, c(t))-t\right].
\end{equation}
\end{defin}

By the triangle inequality, the function $t\rightarrow d(p,
c(t))-t$ is monotonically decreasing.\footnote{In general, not
strictly.} Also by the triangle inequality, it satisfies
$-d(c(0),p)\le d(p, c(t))-t.$ Thus, the limit in
equation~\eqref{busem_eq} exists.

\medskip\noindent
\begin{lem}  \label{busem_lip_lem}
Let $c:\R_+\rightarrow M$ be a ray, and let $p, q\in M$ be
arbitrary points. Then
$$
|B_c(p)-B_c(q)|\le d(p, q).
$$
\begin{proof}
Apply the triangle inequality to the triangle with corners $p,
q,c(t)$, and take the limit $t\to\infty$.
\end{proof}
\end{lem}

By Lemma~\ref{busem_lip_lem}, any Busemann function is lipschitz,
with the lipschitz constant $1$.

\medskip\noindent
\begin{prop}   \label{coray_thm}
Let $c:\R_+\to M$ be a ray. A geodesic $\tc:\R_+\to M$ is a coray
to $c$ if and only if for all $s,t\in\R_+$ the equation
$$
B_c(\tc(t))-B_c(\tc(s))=s-t
$$
holds.
\begin{proof}
This follows from equations (22.16) and (22.20) in \cite{Bu}.
\end{proof}
\end{prop}

\medskip
We use Proposition~\ref{coray_thm} to relax the requirements in
Definition~\ref{ray_def}.

\medskip
%
\begin{lem}  \label{coray_lem}
Let $c:\R_+\to M$ be a ray, and let $C\subset M$ be its trace. Let
$L_n\to\infty$ be a positive sequence. Let $c_n:[0,L_n]\to M$ be
minimal geodesics such that $\lim_{n\to\infty} d(c_n(L_n),C)=0$.

If $\tc:\R_+\to M$ is a geodesic such that $\dot{\tc}(0)$ is a
point of accumulation of the sequence $\dot{c}_n(0)$, then $\tc$
is a coray to $c$.
\begin{proof}
We assume without loss of generality that $\lim_{n\to\infty}
\dot{c}_n(0)=\dot{\tilde{c}}(0)$. For every $n\in\N$ there is a
number $t_n\in[0,\infty)$ such that the sequence
$\ve_n=d(c_n(L_n),c(t_n))$ converges to zero. Since the geodesics
$c_n$ are minimal, the condition $\lim_{n\to\infty}L_n=\infty$
implies that $\lim_{n\to\infty}t_n=\infty$. By the triangle
inequality, for all $n$ and any $t\in[0,L_n]$ we have
\begin{equation}                  \label{coray_eq}
L_n - t - \ve_n\le d(c_n(t), c(t_n))\le L_n - t + \ve_n.
\end{equation}
Let $s, t>0$ be arbitrary. Using that $\lim\dot{c}_n(0) =
\dot{\tilde{c}}(0)$ and equation~\eqref{coray_eq}, we have
\begin{eqnarray*}        
B_c(\tilde{c}(t)) - B_c(\tilde{c}(s))   =
\lim_{n\to\infty} [d(\tilde{c}(t), c(t_n)) - d(\tilde{c}(s), c(t_n))]\\
=\lim_{n\to\infty} [d(c_n(t), c(t_n)) - d(c_n(s),
c(t_n))]\ge\lim_{n\to\infty} (s-t-2\ve_n)=s-t.
\end{eqnarray*}
Combining this inequality with Lemma~\ref{busem_lip_lem}, we
obtain
\begin{equation}        \label{coray_EQ}
B_c(\tilde{c}(t)) - B_c(\tilde{c}(s))=s-t.
\end{equation}
The claim now follows from Proposition~\ref{coray_thm}.
\end{proof}
\end{lem}
\section{Outline of the proof}   \label{idea}
For the benefit of the reader, we will outline the main ideas in
the proof of Theorem~\ref{main_thm}. Let $(T^2,\overline{g})$ be a non-flat
two-dimensional torus; our goal is to find a pair of
points in $(T^2,\overline{g})$ that cannot be blocked
away from each other by a finite blocking set.


\medskip

By a classical theorem of E.~Hopf \cite{Ho}, a riemannian two-torus is flat if and only if
it has no conjugate points. Thus, the
torus $(T^2, \overline{g})$ has conjugate points. Then, by a
theorem of N.~Innami, there exists
a nontrivial free homotopy class $\alpha$ of closed curves such that
$(T^2, \overline{g})$ cannot be foliated by geodesics in  $\alpha$.
See \cite{In}, Corollary 3.2; see also the proof of Theorem 6.1 in \cite{Ba2}.

\medskip

Let $\mm_{\al}$ be the set of geodesics of minimal length in the
class $\alpha$. By results that go back to M.~Morse \cite{Mo} and
G.~Hedlund \cite{He},\footnote{These results provide an important
part of the Aubry-Mather theory \cite{Ba1}.} these geodesics do
not self-intersect and are pairwise disjoint. (Generically, $\mm_{\al}$ consists of a single
geodesic.)

\medskip


The geodesics in $\mm_{\al}$ foliate a compact, proper subset,
$N\subset T^2$.  Let $Z\subset T^2$ be a connected component of
$T^2\setminus N$; let $p,q\in Z$ be any pair of points. We will
show that the pair $p,q$ is insecure, i. e., that we cannot block
$p$ away from $q$ by a finite blocking set.

\medskip

We denote by $(\RR,g)$ the riemannian universal covering; let
$\pi:(\RR,g)\to(T^2,\overline{g})$ be the projection. 
Let $S\subset\RR$ be a connected component of
$\pi^{-1}(Z)$. Then $S$ is an open strip. The boundary $\bo S$ is a
disjoint union of traces of two minimal geodesics, $c_0:\R\to(\RR,g)$ and $c_1:\R\to(\RR,g)$.
Let $C_0,C_1$ be the respective traces; then $\bo S=C_0\cup C_1$.

\medskip

Let $P,Q_0\in S$ be arbitrary points such that $\pi(P)=p,\pi(Q_0)=q$.
Using the action of the stabilizer of $S$ in $\pi_1(T^2)=\ZZ$, we produce an infinite sequence of points 
$Q_1,\dots,Q_n,\ldots\in S$ such that  $\pi(Q_n)=q$ and the sequence of distances $L_n=d(P,Q_n)$
goes to infinity.
Let now $\tc_n:[0,L_n]\to S$ be a sequence of minimal geodesics such that
$\tc_n(0)=P$ and $\tc_n(L_n)=Q_n$. 

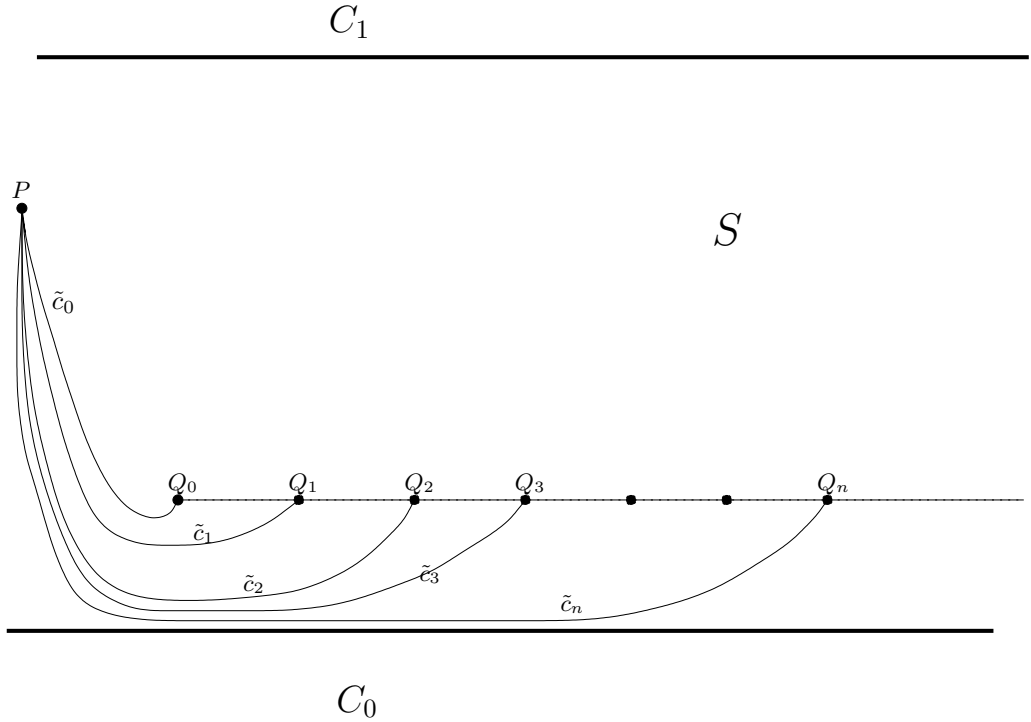
\begin{figure}[htbp]
\begin{center}
\input{toryura.pstex_t}
\caption{A sequence of minimal geodesics in the universal covering whose projections
to the torus cannot be blocked by a finite point set.
}
\label{fig1}
\end{center}
\end{figure}

Lemma~\ref{key_lem} in section~\ref{key}
implies that most of the
time the geodesics $\tc_n$ are close to $\bo S$. More precisely, 
for any $\ve>0$ there exists $T=T(\ve)>0$ such that 
for all $t\in[T,L_n-T]$ the points $\tc_n(t)$ are $\ve$-close to $\bo S$.
Figure~\ref{fig1} illustrates the behavior of this sequence of geodesics.

\medskip

Set $c_n=\pi\circ\tc_n$. Then the geodesics $c_n:[0,L_n]\to Z$
join the points $p,q$. Let $z\in Z\setminus\{p,q\}$ be an arbitrary point.
The preceding discussion implies that at most a finite
number of the geodesics $c_n$ passes through $z$. On the other hand, if $z\in T^2\setminus Z$,
no geodesic in our sequence passes through it.
Thus, any point $z\in T^2$
can block at most a finite number of joining geodesics
in the infinite sequence $c_n$. 
Hence, we cannot block the points $p,q$ away from each other by
a finite set of blocking points.

We will now illustrate the preceding discussion with the example of tori
of revolution. 
\begin{exa}   \label{tor_rev_exa}
{\em  
Let $0<r<R$, and set $C=C(r,R)=\{(x,0,z):(x-R)^2+z^2=r^2\}$.
This is a circle of radius $r$ in the $xz$-plane. The torus of revolution, $T(r,R)\subset\R^3$,
is obtained by revolving $C$ about the $z$-axis. The circles in $T(r,R)$ obtained 
by revolving points in $C$ are the {\em circles of latitude}.

Exactly two of the circles of latitude are geodesics: The inner and the outer equators.
The inner (resp. outer) equator $E_{\inn}$ (resp. $E_{\out}$) is the 
circle of latitude corresponding to the point
$(R-r,0,0)\in C$ (resp. $(R+r,0,0)\in C$). Their lengths are $2\pi(R-r)$ and $2\pi(R+r)$ respectively.
The two equators are freely homotopic; let $\al$ be their homotopy class. 

Thus,  $N=E_{\inn}$, and $Z=T(r,R)\setminus E_{\inn}$.
Note that the set $\mm_{\al}$ consists of a single geodesic; although
the tori of revolution are very special, this is the generic situation for  riemannian tori. 
By the preceding argument, any points
$p,q\in T(r,R)\setminus E_{\inn}$ cannot be blocked away from each other by a finite
blocking set.

}
\end{exa}

\vspace{4mm}

\section{Minimal geodesics in the strip $S$}   \label{constr}
We will use the notation of section~\ref{idea}; in particular, we
use the identification $(T^2, \overline{g})=(\RR,g)/\ZZ$. 
If $S\subset\RR$, we denote
by
$$\stab(S) = \{j\in\ZZ:\,S+j=S\}
$$
the stabilizer of $S$.
Recall that a nonzero vector $j\in\ZZ$ is {\em prime} if there do
not exist $n\in\N$, $n\ge 2$, and $k\in\ZZ$ such that $j=nk$.

%
\begin{prop}     \label{strip_prop}
Let $(T^2, \overline{g})$ be a nonflat riemannian torus. Let
$(\RR, g)$ be its universal covering, and let $\pi:(\RR,
g)\to(T^2, \overline{g})$ be the projection.

Then there exists a connected open set $S\subset\RR$ with 
totally geodesic boundary, such that the following statements
hold.
\begin{itemize}
\item[(a)]
The group $\stab(S)$ is generated by a prime vector.
\item[(b)]
If $j\in\ZZ\setminus\stab(S)$, then
$(S+j)\cap\overline{S}=\emptyset$.
\item[(c)]
The boundary of $S$ has two connected components, say $C_0$ and
$C_1$. There are minimal geodesics $c_0, c_1:\R\to(\RR,g)$ whose
traces are $C_0$ and $C_1$, respectively.
\item[(d)]
Let $c:\R\to S$ be a geodesic such that $\pi\circ c:\R\to T^2$ is
periodic. Then $c$ is not minimal.
\end{itemize}
\begin{proof}
By Corollary 3.2 in \cite{In}, there exists a nontrivial free
homotopy class, say $\alpha$, of closed curves in $T^2$ having the
following property: There does not exist a family of closed
geodesics in the class $\alpha$ whose traces foliate $T^2$. We can
assume that $\alpha$ is prime.



Let $L$ be the minimal length of a curve in $\al$; we denote by
$\mm_{\al}$ the set of closed geodesics in the class $\al$ having
length $L$. Clearly, $\mm_{\al}\ne\emptyset$. By Theorem 6.5 and
Theorem 6.6 in \cite{Ba1}, the trace of every $c\in\mm_{\al}$ is
an embedded curve in $T^2$. Moreover, if $c,\tc\in\mm_{\al}$, then
either their traces are disjoint or $c$ and $\tc$ coincide up to a
translation of the parameter.

\medskip

Let $N$ be the union of the traces of geodesics in $\mm_{\al}$. By
our choice of $\al$, the set $N\subset T^2$ is a proper, nonempty, closed
subset. Let $Z$ be a connected component of $T^2\setminus N$. Let
$\bo Z=\overline{Z}\setminus Z$ be its boundary. Then either $\bo
Z$ is the trace of a geodesic in $\mm_{\al}$ or $\bo Z$ is the
union of traces of two geodesics in  $\mm_{\al}$.\footnote{We
point out that $Z$ is homeomorphic to the cylinder
$S^1\times(0,1)$.}

Let $S$ be a connected component of $\pi^{-1}(Z)\subset\RR$. Then
the boundary of $S$ is the union of the traces of two geodesics
$c_0,c_1:\R\to\RR$ such that $\pi\circ c_0$ and $\pi\circ c_1$
belong to $\mm_{\al}$. Let $k\in\ZZ$ correspond to $\al$. Then for
all $t\in\R$ we have
\begin{equation}   \label{period_eq}
c_0(t+L)=c_0(t)+k,\,c_1(t+L)=c_1(t)+k.
\end{equation}
%

Theorem 6.6 in \cite{Ba1} implies that $c_0$ and $c_1$ are minimal
geodesics. The remaining statements in (a), (b), and (c) now follow by
elementary topological arguments; claim (d) follows from Theorem
6.7 in \cite{Ba1}.
\end{proof}
\end{prop}

\medskip

\section{The key lemma}   \label{key}
We will use the setting and the notation of
Proposition~\ref{strip_prop}. The following statement is
crucial in our proof of Theorem~\ref{main_thm}. 
We will refer to it as the Key Lemma.

\begin{lem}     \label{key_lem}
For any $\ve>0$ and any $\delta>0$ there exists
$T=T(\ve,\delta)>0$ such that the following holds. If
$c:[0,L]\rightarrow S$ is a minimal geodesic and $d(c(0),\partial
S)\ge\delta$ then $d(c(t),\partial S)\le\ve$ for all
$t\in[T,L-T]$.
\end{lem}


The proof of Lemma~\ref{key_lem} is based on the results of
M.~Morse \cite{Mo} about minimal geodesics in $S$ and on a result
from \cite{Ba2} concerning the rays in $S$.
We need a few technical lemmas.

\begin{lem}     \label{dist_lem}
Let  $c:[0,\infty)\rightarrow S$ be a ray. Then $\lim_{t\to\infty}
d(c(t),\bo S)=0$.
\begin{proof}
The claim follows from Theorem 3.7 in \cite{Ba2}, interpreted as a
statement about minimal geodesics in $(\RR,g)$. See Example (1) on page 51
in \cite{Ba2} for details.
\end{proof}
\end{lem}
\medskip

Throughout this section we will use the following notational conventions.
With any geodesic $c:\R\to\overline{S}$
we will associate two geodesics $c_{\pm}:\R_+\to\overline{S}$ as follows.
The geodesic $c_+$ is the restriction of $c$ to the positive half-line.
We define the geodesic $c_-$ by $c_-(t)=c(-t)$.


We will denote by $C,C_+,C_-\subset\overline{S}$ the respective traces of
$c,c_+,c_-$.

\vspace{6mm}

\begin{defin}     \label{coherent_def}
Let $c_0,c_1:\R\to\overline{S}$ be two geodesics such that their
traces $C_0,C_1$ are the two components of $\bo S$. We say that
the geodesics $c_0,c_1$ are {\em coherently oriented} if for any
time sequence $t_n\to\infty$ the two point sequences
$c_0(t_n),c_1(t_n)\in\bo S$ converge to the same end of
$\overline{S}$.
\end{defin}

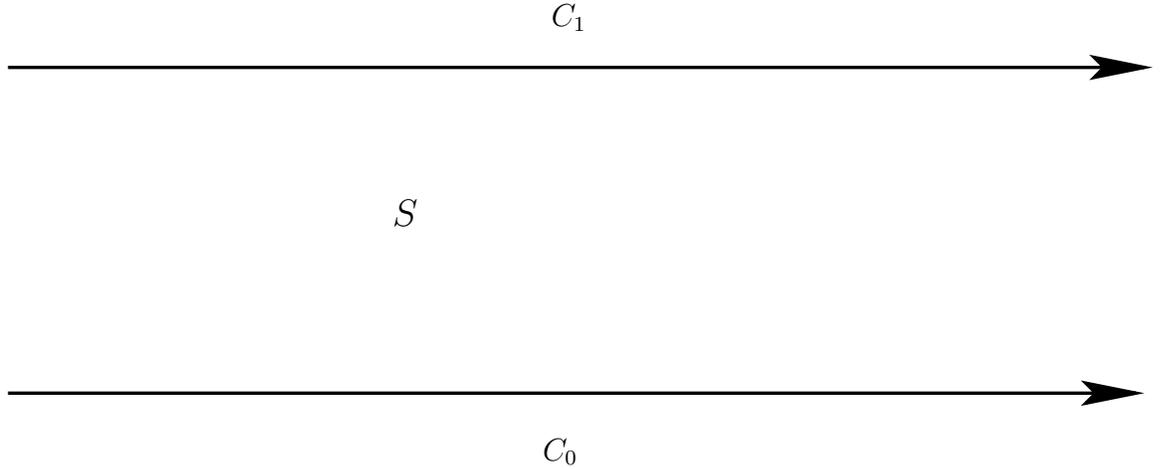
\begin{figure}[htbp]
\begin{center}
\input{tori2.pstex_t}
\caption{A strip with coherently oriented boundary components.}
\label{fig2}
\end{center}
\end{figure}

Figure~\ref{fig2} illustrates Definition~\ref{coherent_def}.
We will also say that $c_0,c_1:\R\to\overline{S}$ are {\em coherent parameterizations}
of $\bo S$.

\medskip

\begin{lem}     \label{component_lem}
Let $\bo S =C_0\cup C_1$, where $c_0,c_1:\R\to\overline{S}$ are
coherently oriented. Let the geodesics $c_{0,\pm}:\R_+\to\bo S$
and $c_{1,\pm}:\R_+\to\bo S$ be as above; let $C_{0,\pm}\subset\bo
S$ and $C_{1,\pm}\subset\bo S$ be the respective traces.

Let now $c:\R\to S$ be a minimal geodesic.
Then, switching $c_0$ with $c_1$ and reversing the orientation of $c$, if need be, we have
\begin{equation}    \label{ends_eq}
\lim_{t\to-\infty}d(c(t),C_{0,-})=0,\ \lim_{t\to+\infty}d(c(t),C_{1,+})=0.
\end{equation}
\begin{proof}
By Lemma~\ref{dist_lem}, $c(t)$ converges to $\bo S$ as $|t|\to\infty$.
By Theorem 15 in \cite{Mo} or by Theorem 6.7 in \cite{Ba1}, the equation
$\lim_{t\to-\infty}d(c(t),C_{0,-})=0$ implies $\lim_{t\to+\infty}d(c(t),C_{1,+})=0$.
\end{proof}
\end{lem}
%
%
%
\begin{rem}    \label{connect_rem}
{\em Let $M$ be a riemannian manifold. If $c:I\to M$ is a geodesic, its {\em inverse} is the
geodesic $c^{-1}:-I\to M$ defined by $c^{-1}(t)=c(-t)$. Lemma~\ref{component_lem} is 
equivalent to the following geometric fact.

\noindent Let $Z\subset T^2$ be as in section~\ref{constr}. Assume, 
for simplicity of exposition, that the closure of $Z$ is a proper subset of $T^2$.
Let $\overline{c}_0,\overline{c}_1:\R\to(T^2,\overline{g})$ be the periodic geodesics 
in the homotopy class $\al$ whose
respective traces are the two components of the boundary $\bo Z$. 

Let now $c:\R\to Z$ be a geodesic whose lift $\tc:\R\to S$ is minimal.
Then  $c$ is  a heteroclinic connection either between 
$\overline{c}_0$ and $\overline{c}_1$ or 
between $\overline{c}_0^{-1}$ and $\overline{c}_1^{-1}$.

}
\end{rem} 

%
%
\medskip

Our next lemma says that if a ray in $S$ is a coray  to a ray
in a boundary component of $S$, then it is asymptotic to this component.

\begin{figure}[htbp]
\begin{center}
\input{tori3.pstex_t}
\caption{Illustration to the proof of Lemma~\ref{asympt_lem}: A configuration that cannot take place.}
\label{fig3}
\end{center}
\end{figure}
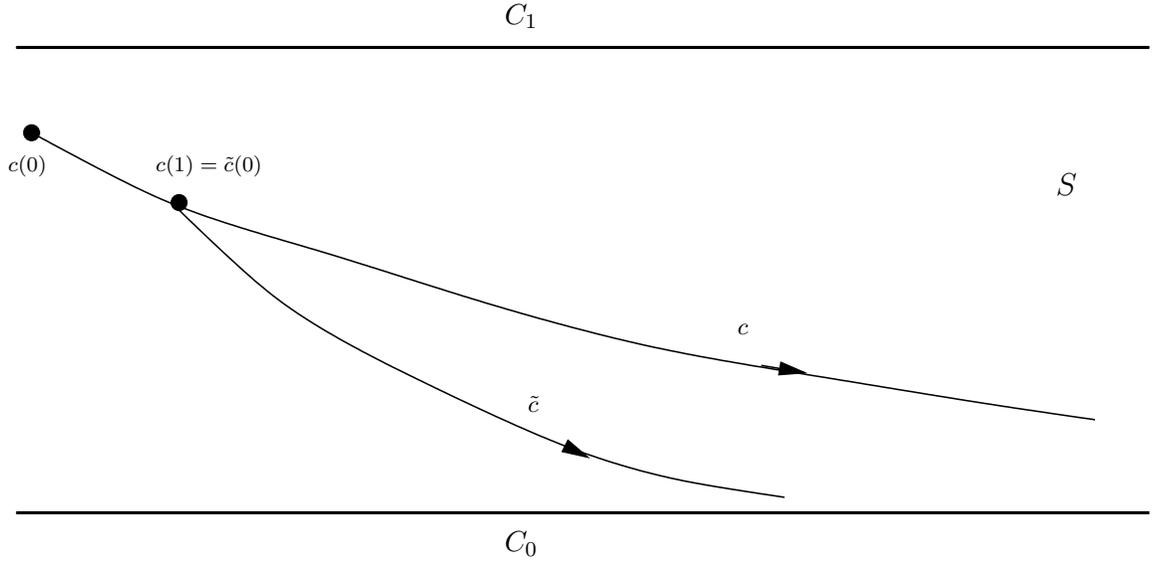

\begin{lem}     \label{asympt_lem}
Let $c_0:\R\to\overline{S}$ be a geodesic whose trace is one of the components of  $\bo S$.

Let $c:\R_+\to S$ be a coray to $c_{0,+}$. Then $\lim_{t\to\infty}d(c(t),C_{0,+})=0$.
\begin{proof}
By Theorem 3.7 in \cite{Ba2}, for any $q\in S$ there exists a ray $\tc:\R_+\to S$ such that
$\tc(0)=q$ and $\lim_{t\to\infty}d(\tc(t),C_{0,+})=0$. Hence, 
by Lemma~\ref{dist_lem} and Lemma~\ref{coray_lem},
$\tc$ is a coray to $c_{0,+}$.

Set $q=c(1)$, and let $\tc:\R_+\to S$ be as above. Thus, both $c$
and $\tc$ are corays to $c_{0,+}$; by construction, $\tc(0)=c(1)$.
The geodesic $t\mapsto c(1+t)$ is also a coray  to $c_{0,+}$ starting at
$q=c(1)=\tc(0)$. By Theorem 22.19 in \cite{Bu} or, by Corollary
3.8 in \cite{Ba2}, there is only one coray to $c_{0,+}$ starting
at $q$. Figure~\ref{fig3} shows a hypothetical configuration of
the rays $c$ and $\tc$ which cannot materialize.

Therefore, the ray $\tc$ satisfies
$\tc(t)=c(1+t)$. Since $\tc$ is asymptotic to $C_{0,+}$, the claim follows.
\end{proof}
\end{lem}
%



\medskip

We will now prove a preliminary variant of the Key Lemma.

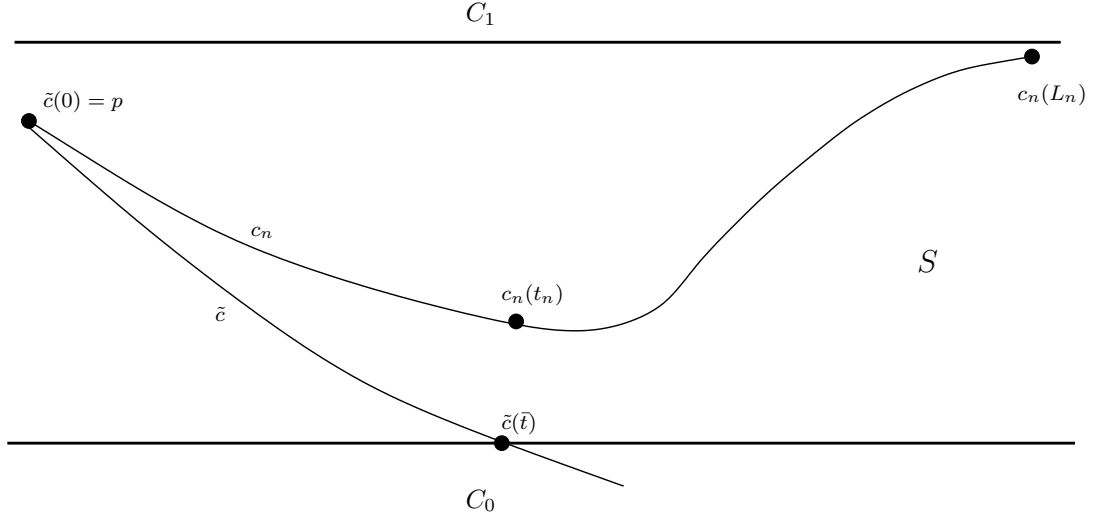
\begin{figure}[htbp]
\begin{center}
\input{tori4.pstex_t}
\caption{Illustration to the proof of Lemma~\ref{prekey_lem}: The behavior of the 
geodesic $\tc$ deduced from the assumption that $\lim_{n\to\infty}t_n<\infty$.}
\label{fig4}
\end{center}
\end{figure}

\begin{lem}     \label{prekey_lem}
For any $\delta>0$ there exists $\eta=\eta(\delta)>0$ such that
the following holds. Let $0<L<\infty$, and let
$c:[0,L]\to S$ be a minimal geodesic such
that $d(c(0),\partial S)\ge\delta$ and $d(c(L),C_1)<\eta$. Then
$d(c(t),C_0)\ge\eta$ for all $t\in[0,L]$.
\begin{proof}
Suppose that the claim fails. Then there exists $\delta>0$, a
sequence of minimal geodesics $c_n:[0,L_n]\to S$, and a sequence
$t_n\in[0,L_n]$ such that $d(c_n(0),\partial S)\ge\delta$,
$\lim_{n\to\infty}d(c_n(L_n),C_1)=0$, and
$\lim_{n\to\infty}d(c_n(t_n),C_0)=0$.

The closed strip $\overline{S}$ is invariant under the group $\stab(S)\simeq k\Z$
that acts on $\overline{S}$ by isometries. We have denoted this action by $z\mapsto z+rk$.
We will use the same notation for the corresponding action of  $\stab(S)$ on geodesics
in $S$. Then, for any
integers $r_1,\dots,r_n,\ldots\in\Z$ the sequence of geodesics $\tc_n=c_n+r_nk$ satisfies
the above conditions. In view of this observation, and the compactness of the
quotient $\overline{S}/k\Z$, we assume without loss of generality
that the vectors $\dot{c}_n(0)$ converge to a limit vector, $v\in T^1(S,g)$;
let $p\in S$ be its footpoint.

We will now prove that $\lim_{n\to\infty}t_n=\infty$. If this
fails, then, by passing to a subsequence, if need be, we have
$\lim t_n=\overline{t}<\infty$. Let $\tc:\R\to\RR$ be the geodesic
with the initial vector $v$. Then $\tc(0)=p\in S$, and
$\tc(\overline{t})=q=\lim_{n\to\infty} c_n(t_n)\in C_0\subset\bo S$. 
Since $\bo S$ is geodesic,
$\tc$ intersects it transversally at $q$. Thus, for $t>\overline{t}$ and
sufficiently close to $\overline{t}$, we have $\tc(t)\notin\overline{S}$.
Figure~\ref{fig4} illustrates the analysis.

On the other hand, $\liminf L_n\ge\overline{t}+d(C_0,C_1)$ implies that
$\tc(t)\in\overline{S}$ for all $t\in[0,\overline{t}+d(C_0,C_1)]$. In view of this
contradiction, $\lim t_n=\infty$.

\medskip

By Lemma~\ref{coray_lem}, the relationships $\lim d(c_n(t_n),C_0)=0$ and
$\lim t_n=\infty$ together imply that $\tc$ is a coray to
$c_{0,+}$ or $c_{0,-}$. Similarly, the relationships $\lim d(c_n(L_n),C_1)=0$ and
$\lim L_n=\infty$ imply that $\tc$ is a coray to $c_{1,+}$ or $c_{1,-}$.
By Lemma~\ref{asympt_lem}, this is impossible.
\end{proof}
\end{lem}

\vspace{3mm}

We will now prove the Key Lemma. Recall that we view geodesics in
$\overline{S}$ as mappings $c:I\to\overline{S}$ of nontrivial
intervals $I\subset\R$. For $t\in I$ the velocity vectors
$\dot{c}(t)$ are unit tangent vectors in $(\RR,g)$. Thus,
$\len(c)=|I|$. If $0\in I$, we
will refer to $\dot{c}(0)$ as the {\em initial vector} of $c$.

\medskip
\begin{proof}({of Lemma~\ref{key_lem}})
Assume that the claim fails. Then for some $\ve>0$, $\delta>0$
there exists  a sequence of minimal geodesics $c_n:[0,L_n]\to S$ such that
the following conditions are satisfied:

\noindent i) For all $n\in\N$ we have $d(c_n(0),\bo S)\ge\delta$;

\noindent ii) For each $n$ there
is $t_n\in[n,L_n-n]$ so that $d(c_n(t_n),\bo S)>\ve$.

As in the proof of Lemma~\ref{prekey_lem}, we assume
without loss of generality that the velocity vectors $\dot{c}_n(t_n)$
converge to a vector $v\in T^1(\overline{S},g)$. Let $\tc:\R\rightarrow S$
be the geodesic such that $v=\dot{\tc}(0)$. Since all of $c_n:[0,L_n]\to S$
are minimal, and $t_n\in[n,L_n-n]$, we conclude that $\tc:\R\rightarrow S$
is a minimal geodesic. By construction, it satisfies
$d(\tc(0),\bo S)\ge\delta$.

Let $\eta=\eta(\delta)>0$ be as in
Lemma~\ref{prekey_lem}. By Lemma~\ref{component_lem}, there are
$s_0, s_1\in\R$ such that $d(\tc(s_0),C_0)<\eta$ and
$d(\tc(s_1),C_1)<\eta$.

Interchanging $C_0$ and $C_1$, if need be, we may assume that
$s_0<s_1$.
For any $t\in\R$ we have $\lim_{n\to\infty}c_n(t_n+t)=\tc(t)$.
In particular, $\tc(s_0)=\lim_{n\to\infty}c_n(t_n+s_0)$ and
$\tc(s_1)=\lim_{n\to\infty}c_n(t_n+s_1)$.
Therefore, for sufficiently large $n$  the inequalities $d(c_n(t_n+s_0), C_0)<\eta$
and $d(c_n(t_n+s_1), C_1)<\eta$ hold. Besides, for sufficiently large $n$  we have $0<t_n+s_0<t_n+s_1<L_n$.

Let $n\in\N$ be any index such that the above conditions hold. Set $L=t_n+s_1$, and let
$c:[0,L]\to S$ be the restriction of $c_n$ to $[0,t_n+s_1]$. Then $d(c(0),\bo S)>\delta$ and
$d(c(L),C_1)<\eta$. But we also have $d(c(t_n+s_0),C_0)<\eta$, and $t_n+s_0\in(0,L)$.
By Lemma~\ref{prekey_lem}, this is impossible.
\end{proof}
\section{Nonflat two-tori are insecure}  \label{proof}
We use the setting of section~\ref{constr} and the notation of section~\ref{key}.
First, we need a technical lemma.
\begin{lem}    \label{vectors_lem}
Let $p,q\in Z$ be arbitrary points. Let $P,Q\in S$ be such
that $\pi(P)=p$, $\pi(Q)=q$. For $n\in\N$ let
$\tc_n:[0,L_n]\to S$ be a minimal geodesic such that
$\tc_n(0)=P$ and $\tc_n(L_n)=Q+nk$.

Then in the sequence of unit tangent vectors $\dot{\tc}_n(0)\in
T_{P}^1(S,g)$ every vector occurs at most a finite number of times.
\begin{proof}
Assume the opposite. Then, by passing to a subsequence of indices,
if need be, we find a unit vector $v\in T_{P}^1(S,g)$ and a sequence
of minimal geodesics $\hc_i:[0,T_i]\to S$ such that
$\dhc_i(0)=v,\hc_i(T_i)=Q+n(i)k$.

We have $\lim_{i\to\infty}n(i)=\infty$. Thus, $\lim_{i\to\infty}\len(\hc_i)=\infty$.
Since all of $\hc_i$ have the same initial vector, the geodesic $\hc_{i+1}$
extends $\hc_i$ from $[0,T_i]$ to $[0,T_{i+1}]$. Therefore, the limit geodesic
$\lim_{i\to\infty}\hc_i=\hc:[0,\infty)\to S$ coincides with $\hc_i$ on $[0,T_i]$.
Hence, $\hc:[0,\infty)\to S$ is a ray.

By Lemma~\ref{dist_lem}, $\hc$ is asymptotic to $\bo S$ at infinity.
On the other hand, for $i\in\N$ we have
$$
d(\hc(T_i),\bo S)=d(Q+n(i)k,\bo S)=d(Q,\bo S)>0.
$$
We have arrived at a contradiction.
\end{proof}
\end{lem}

\medskip

Theorem~\ref{main_thm} will follow immediately from the proposition below.


\begin{prop}       \label{no_block_prop}
Let $p,q\in Z$ be arbitrary points.
Then they cannot be blocked away from each other by a finite blocking set.
\begin{proof}
Let $P,Q\in S$ be such that $\pi(P)=p$,
$\pi(Q)=q$. For $n\in\N$ let $\tc_n:[0,L_n]\to S$ be a minimal geodesic
such that $\tc_n(0)=P$ and $\tc_n(L_n)=Q+nk$.
Set $c_n=\pi\circ\tc_n$. We will show that no point
belongs to the interior of  infinitely many geodesics $c_n$.

\medskip

Suppose this is false. Then  there is a point $z\in
T^2\setminus\{p,q\}$, an infinite set $I\subset\N$, and a function
$i\mapsto t_i\in(0,L_i)$ on $I$ such that $c_i(t_i)=z$. We will
now analyze all apriori possible behaviors of the sequence $t_i$,
as $i\to\infty$.

\medskip

Since all of the geodesics in question belong to $Z$, we have
$z\in Z\setminus\{p,q\}$. By construction, $d(p,z)\le t_i$
and $d(z,q)\le L_i-t_i$. Thus the sequences $t_i$ and $L_i-t_i$
are bounded away from zero.

\medskip

Suppose first that $\limsup_{i\to\infty}t_i<\infty$. Then, by passing to an appropriate subsequence of $I$,
if need be, we obtain the following situation: $\lim_{i\to\infty}t_i=T\in(0,\infty)$ and the
vectors $\dc_i(0)$ converge. Let $w=\lim_{i\to\infty}\dc_i(0)$. Let $c:\R_+\to Z$ be
the geodesic such that $c(0)=p$ and $\dc(0)=w$. Let $\tc:\R_+\to S$ be its lift such that
$\tc(0)=P$. Since $\tc$ is a limit of minimal geodesics, it is a ray.
Thus, $c$ has no conjugate points.
On the other hand, $c:[0,T]\to Z$ is the limit of the geodesics $c_i:[0,t_i]\to Z$.
We have $c_i(0)=p=c(0)$ and $c_i(t_i)=z=c(T)$. By Lemma~\ref{vectors_lem}, we can assume that
the vectors $\dc_i(0)$ are distinct.
Thus, the points $p$ and $z$ are conjugate along $c$; we have arrived at a contradiction.

We reduce the  case $\limsup_{i\to\infty}(L_i-t_i)<\infty$ to the
preceding one by switching the roles of the points $p$ and $q$ and
simultaneously reversing the directions of the geodesics $c_i$. We
conclude that $\limsup_{i\to\infty}(L_i-t_i)<\infty$ is
impossible as well.

\medskip

The only remaining possibility is $\lim_{i\to\infty}t_i =
\lim_{i\to\infty}(L_i - t_i)=\infty$. Then, by
Lemma~\ref{key_lem}, $\lim_{i\to\infty}d(\tc_i(t_i),\bo S)=0$.
Since $z=\pi(\tc_i(t_i))$ for all $i$, we conclude that
$z\in\pi(\bo S)\cap\pi(S)=\pi(\bo S)\cap Z$. This contradicts claim (b) in
Proposition~\ref{strip_prop}.

\medskip

We have examined all possibilities for the sequence
$t_i\in(0,L_i),\,i\in I,$ and
arrived at a contradiction in each case. Therefore, at most a finite number of the
geodesics $c_n$ pass through any point in $T^2\setminus\{p,q\}$. 
\end{proof}
\end{prop}

\medskip

\noindent{\em Proof of Theorem~\ref{main_thm}.} By \cite{Gut05} or \cite{GS06},
a flat torus is secure. By Proposition~\ref{no_block_prop}, a nonflat two-torus
is insecure.\qed

\medskip

We point out that the flat tori are distinguished amongst all riemannian two-tori
by the security of pairs $y=x$.
\begin{corol}    \label{flat_cor1}
A two-dimensional riemannian torus is flat if and only if all pairs $x,x$ are secure.
\begin{proof}
It suffices to show that a nonflat riemannian two-torus contains at least one point
that cannot be blocked away from itself. Let $Z$ be the cylinder from Proposition~\ref{no_block_prop}.
Then any point $x\in Z$ cannot be blocked away from itself.
\end{proof}
\end{corol}

\medskip


%
%

\end {document}

%% file: toryura.pstex_t
\begin{picture}(0,0)%
\includegraphics{toryura.pstex}%
\end{picture}%
\setlength{\unitlength}{3355sp}%
\begingroup\makeatletter\ifx\SetFigFont\undefined%
\gdef\SetFigFont#1#2#3#4#5{%
  \reset@font\fontsize{#1}{#2pt}%
  \fontfamily{#3}\fontseries{#4}\fontshape{#5}%
  \selectfont}%
\fi\endgroup%
\begin{picture}(7618,5257)(5,-5449)
\put( 75,-1588){\makebox(0,0)[lb]{\smash{{\SetFigFont{9}{10.8}{\rmdefault}{\mddefault}{\updefault}{\color[rgb]{0,0,0}$P$}%
}}}}
\put(1229,-3746){\makebox(0,0)[lb]{\smash{{\SetFigFont{9}{10.8}{\familydefault}{\mddefault}{\updefault}{\color[rgb]{0,0,0}$Q_0$}%
}}}}
\put(2121,-3746){\makebox(0,0)[lb]{\smash{{\SetFigFont{9}{10.8}{\familydefault}{\mddefault}{\updefault}{\color[rgb]{0,0,0}$Q_1$}%
}}}}
\put(2977,-3746){\makebox(0,0)[lb]{\smash{{\SetFigFont{9}{10.8}{\familydefault}{\mddefault}{\updefault}{\color[rgb]{0,0,0}$Q_2$}%
}}}}
\put(3795,-3746){\makebox(0,0)[lb]{\smash{{\SetFigFont{9}{10.8}{\familydefault}{\mddefault}{\updefault}{\color[rgb]{0,0,0}$Q_3$}%
}}}}
\put(6027,-3746){\makebox(0,0)[lb]{\smash{{\SetFigFont{9}{10.8}{\familydefault}{\mddefault}{\updefault}{\color[rgb]{0,0,0}$Q_n$}%
}}}}
\put(2419,-360){\makebox(0,0)[lb]{\smash{{\SetFigFont{14}{16.8}{\familydefault}{\mddefault}{\updefault}{\color[rgb]{0,0,0}$C_1$}%
}}}}
\put(5251,-1936){\makebox(0,0)[lb]{\smash{{\SetFigFont{17}{20.4}{\rmdefault}{\mddefault}{\updefault}{\color[rgb]{0,0,0}$S$}%
}}}}
\put(2476,-5386){\makebox(0,0)[lb]{\smash{{\SetFigFont{14}{16.8}{\familydefault}{\mddefault}{\updefault}{\color[rgb]{0,0,0}$C_0$}%
}}}}
\put(373,-2407){\makebox(0,0)[lb]{\smash{{\SetFigFont{10}{12.0}{\familydefault}{\mddefault}{\updefault}{\color[rgb]{0,0,0}$\tc_0$}%
}}}}
\put(1415,-4118){\makebox(0,0)[lb]{\smash{{\SetFigFont{9}{10.8}{\familydefault}{\mddefault}{\updefault}{\color[rgb]{0,0,0}$\tc_1$}%
}}}}
\put(1787,-4490){\makebox(0,0)[lb]{\smash{{\SetFigFont{9}{10.8}{\familydefault}{\mddefault}{\updefault}{\color[rgb]{0,0,0}$\tc_2$}%
}}}}
\put(3089,-4416){\makebox(0,0)[lb]{\smash{{\SetFigFont{9}{10.8}{\familydefault}{\mddefault}{\updefault}{\color[rgb]{0,0,0}$\tc_3$}%
}}}}
\put(4130,-4639){\makebox(0,0)[lb]{\smash{{\SetFigFont{9}{10.8}{\familydefault}{\mddefault}{\updefault}{\color[rgb]{0,0,0}$\tc_n$}%
}}}}
\end{picture}%

%% file: tori2.pstex_t
\begin{picture}(0,0)%
\includegraphics{tori2.pstex}%
\end{picture}%
\setlength{\unitlength}{2763sp}%
\begingroup\makeatletter\ifx\SetFigFont\undefined%
\gdef\SetFigFont#1#2#3#4#5{%
  \reset@font\fontsize{#1}{#2pt}%
  \fontfamily{#3}\fontseries{#4}\fontshape{#5}%
  \selectfont}%
\fi\endgroup%
\begin{picture}(10341,4131)(868,-6499)
\put(5776,-2536){\makebox(0,0)[lb]{\smash{{\SetFigFont{12}{14.4}{\familydefault}{\mddefault}{\updefault}{\color[rgb]{0,0,0}$C_1$}%
}}}}
\put(5701,-6436){\makebox(0,0)[lb]{\smash{{\SetFigFont{12}{14.4}{\familydefault}{\mddefault}{\updefault}{\color[rgb]{0,0,0}$C_0$}%
}}}}
\put(4351,-4336){\makebox(0,0)[lb]{\smash{{\SetFigFont{14}{16.8}{\rmdefault}{\mddefault}{\updefault}{\color[rgb]{0,0,0}$S$}%
}}}}
\end{picture}%

%% file: tori3.pstex_t
\begin{picture}(0,0)%
\includegraphics{tori3.pstex}%
\end{picture}%
\setlength{\unitlength}{2565sp}%
\begingroup\makeatletter\ifx\SetFigFont\undefined%
\gdef\SetFigFont#1#2#3#4#5{%
  \reset@font\fontsize{#1}{#2pt}%
  \fontfamily{#3}\fontseries{#4}\fontshape{#5}%
  \selectfont}%
\fi\endgroup%
\begin{picture}(11058,5331)(976,-6949)
\put(976,-3211){\makebox(0,0)[lb]{\smash{{\SetFigFont{8}{9.6}{\familydefault}{\mddefault}{\updefault}{\color[rgb]{0,0,0}$c(0)$}%
}}}}
\put(5776,-1786){\makebox(0,0)[lb]{\smash{{\SetFigFont{11}{13.2}{\familydefault}{\mddefault}{\updefault}{\color[rgb]{0,0,0}$C_1$}%
}}}}
\put(8026,-4786){\makebox(0,0)[lb]{\smash{{\SetFigFont{9}{10.8}{\familydefault}{\mddefault}{\updefault}{\color[rgb]{0,0,0}$c$}%
}}}}
\put(5776,-6886){\makebox(0,0)[lb]{\smash{{\SetFigFont{11}{13.2}{\familydefault}{\mddefault}{\updefault}{\color[rgb]{0,0,0}$C_0$}%
}}}}
\put(6001,-5536){\makebox(0,0)[lb]{\smash{{\SetFigFont{9}{10.8}{\familydefault}{\mddefault}{\updefault}{\color[rgb]{0,0,0}$\tc$}%
}}}}
\put(2401,-3211){\makebox(0,0)[lb]{\smash{{\SetFigFont{8}{9.6}{\familydefault}{\mddefault}{\updefault}{\color[rgb]{0,0,0}$c(1)=\tc(0)$}%
}}}}
\put(11101,-3436){\makebox(0,0)[lb]{\smash{{\SetFigFont{12}{14.4}{\rmdefault}{\mddefault}{\updefault}{\color[rgb]{0,0,0}$S$}%
}}}}
\end{picture}%

%% file: tori4.pstex_t
\begin{picture}(0,0)%
\includegraphics{tori4.pstex}%
\end{picture}%
\setlength{\unitlength}{2368sp}%
\begingroup\makeatletter\ifx\SetFigFont\undefined%
\gdef\SetFigFont#1#2#3#4#5{%
  \reset@font\fontsize{#1}{#2pt}%
  \fontfamily{#3}\fontseries{#4}\fontshape{#5}%
  \selectfont}%
\fi\endgroup%
\begin{picture}(11654,5331)(943,-6949)
\put(5776,-1786){\makebox(0,0)[lb]{\smash{{\SetFigFont{10}{12.0}{\familydefault}{\mddefault}{\updefault}{\color[rgb]{0,0,0}$C_1$}%
}}}}
\put(3526,-4036){\makebox(0,0)[lb]{\smash{{\SetFigFont{8}{9.6}{\familydefault}{\mddefault}{\updefault}{\color[rgb]{0,0,0}$c_n$}%
}}}}
\put(1351,-2686){\makebox(0,0)[lb]{\smash{{\SetFigFont{8}{9.6}{\familydefault}{\mddefault}{\updefault}{\color[rgb]{0,0,0}$\tc(0)=p$}%
}}}}
\put(6151,-4711){\makebox(0,0)[lb]{\smash{{\SetFigFont{8}{9.6}{\familydefault}{\mddefault}{\updefault}{\color[rgb]{0,0,0}$c_n(t_n)$}%
}}}}
\put(11551,-2611){\makebox(0,0)[lb]{\smash{{\SetFigFont{8}{9.6}{\familydefault}{\mddefault}{\updefault}{\color[rgb]{0,0,0}$c_n(L_n)$}%
}}}}
\put(6151,-6061){\makebox(0,0)[lb]{\smash{{\SetFigFont{8}{9.6}{\familydefault}{\mddefault}{\updefault}{\color[rgb]{0,0,0}$\tc(\bar{t})$}%
}}}}
\put(3151,-4936){\makebox(0,0)[lb]{\smash{{\SetFigFont{8}{9.6}{\rmdefault}{\mddefault}{\updefault}{\color[rgb]{0,0,0}$\tc$}%
}}}}
\put(10501,-4411){\makebox(0,0)[lb]{\smash{{\SetFigFont{12}{14.4}{\rmdefault}{\mddefault}{\updefault}{\color[rgb]{0,0,0}$S$}%
}}}}
\put(5776,-6886){\makebox(0,0)[lb]{\smash{{\SetFigFont{10}{12.0}{\familydefault}{\mddefault}{\updefault}{\color[rgb]{0,0,0}$C_0$}%
}}}}
\end{picture}%